\documentclass[11pt]{amsart}

\usepackage{amsmath}
\usepackage{amssymb}
\usepackage[all]{xy}

\makeatletter
\@addtoreset{equation}{section}
\makeatother

\def\Ree{{\mathbb R}}

\def\fA{{\mathcal A}}

\def\fC{{\mathcal C}}

\def\fI{{\mathcal I}}
\def\fJ{{\mathcal J}}

\def\fS{{\mathcal S}}

\def\fV{{\mathcal V}}
\def\fW{{\mathcal W}}
\def\fX{{\mathcal X}}
\def\fY{{\mathcal Y}}

\def\Del{\Delta}
\def\eps{\varepsilon}
\def\gam{\gamma}
\def\Gam{\Gamma}

\def\sig{\sigma}

\def\vphi{\varphi}

\def\supp{\mathrm{supp}}

\def\comp{\raisebox{.2ex}{${\scriptstyle\circ}$}}
\def\con{\negmedspace\ast\negmedspace}
\def\cross{\negmedspace\times\negmedspace}
\def\mult{\negmedspace\cdot\negmedspace}
\def\setdif{\setminus}
\def\to{\rightarrow}

\def\norm#1{\left\|{#1}\right\|}

\def\unorm#1{\left\|#1\right\|_\infty}
\def\wbar#1{\overline{#1}}
\def\what#1{\widehat{#1}}

\def\aand{\text{ and }}

\def\ffor{\text{ for }}

\def\iin{\text{ in }}
\def\oof{\text{ of }}
\def\oor{\text{ or }}
\def\wwhere{\text{ where }}

\def\endpf{{\hfill$\square$\medskip}}
\def\proof{{\noindent{\bf Proof.}\thickspace}}

\def\bloneg{\mathrm{L}^1(G)}

\def\falg{\mathrm{A}(G)}
\def\falgg{\mathrm{A}(G\cross G)}
\def\fal#1{\mathrm{A}(#1)}

\def\fdelg{\mathrm{A}_{\Del}(G)}

\def\fdel#1{\mathrm{A}_{\Del}(#1)}
\def\fdelgg{\mathrm{A}_{\Del}(G\cross G)}

\def\fdelbg#1{\mathrm{A}_{\Del^{#1}}(G)}

\def\fgamg{\mathrm{A}_{\gam}(G)}

\def\hull{\mathrm{hull}}

\def\ideal{\mathrm{I}}

\def\z#1{\mathrm{Z}(#1)}

\def\fgamnorm#1{\left\|#1\right\|_{\mathrm{A}_{\gam}}}

\begin{document}

\newtheorem{ideals}{Lemma}[section]

\newtheorem{weakamen}{Theorem}[section]
\newtheorem{idealweakamen}[weakamen]{Theorem}
\newtheorem{weakamenfdelg}[weakamen]{Theorem}
\newtheorem{amenfdelg}[weakamen]{Theorem}
\newtheorem{inversion}[weakamen]{Proposition}
\newtheorem{specsynth}[weakamen]{Lemma}
\newtheorem{antidiag}[weakamen]{Theorem}
\newtheorem{johnson}[weakamen]{Proposition}

\newtheorem{genweakamen}{Theorem}[section]
\newtheorem{localproperty}[genweakamen]{Lemma}
\newtheorem{locantidiag}[genweakamen]{Theorem}
\newtheorem{locantidiag1}[genweakamen]{Corollary}

\title[Weak amenability of Fourier algebras]
{Weak amenability of Fourier algebras on compact groups
}

\author{Brian E. Forrest, Ebrahim Samei and Nico Spronk}

\begin{abstract}
We give for a compact group $G$, a full characterisation of when
its Fourier algebra $\falg$
is weakly amenable:  when the connected component
of the identity $G_e$ is abelian.  This condition is also
equivalent to the hyper-Tauberian property for $\falg$,  
and to having the anti-diagonal $\check{\Del}=\{(s,s^{-1}):s\in G\}$
being a set of spectral synthesis for  $\falgg$.
We show the relationship between amenability and weak amenability
of $\falg$, and (operator) amenability and (operator) weak amenability
of $\fdelg$, an algebra defined by the authors in \cite{forrestss}.
We close by extending our results to some classes of non-compact, locally
compact groups, including small invariant neighbourhood groups and 
maximally weakly almost periodic groups.
\end{abstract}

\maketitle

\footnote{{\it Date}: \today.


Research of the the first named author supported by NSERC Grant 90749-00.
Research of the second named author supported by an NSERC Post Doctoral 
Fellowship.  
Research of the the third named author supported by NSERC Grant 312515-05.}

\section{Background and Notation}

\subsection{Context}
B. E. Johnson \cite{johnson}, showed that for the compact Lie group
$G=\mathrm{SO}(3)$ , the Fourier algebra $\falg$ is not weakly
amenable.  At the time that article was written, the result was a surprise,
since it was expected that $\falg$ should be amenable for any amenable locally
compact group, and in particular for any compact group, in analogy to the famous
characterisation of amenability of group algebras $\bloneg$ \cite{johnsonM}.  
Very soon after that, Z.-J. Ruan \cite{ruan} showed that
if one puts a suitable operator space structure on $\falg$, then
$\falg$ is operator amenable exactly when $G$ is amenable. This was the first evidence 
of the important role that the operator space structure would play in 
understanding the cohomology of $\falg$ for non-commuative groups. Indeed, Johnson's paper 
combined with earlier work of V. Losert \cite{losert}, led to the conjecture that $\falg$ 
would be amenable as a Banach algebra precisely when $G$ was {\it almost abelian},
i.e.\ $G$ admits an abelian subgroup of finite index. 
This conjecture was finally settled affirmatively by the first author and V. Runde
\cite{forrestr}.  Refelected in this result is the fact that the Fourier
algebra has a trivial operator space structure precisely when $G$ admits an
abelian subgroup of finite index \cite{forrestw}.

On the separate question of weak amenability, Johnson \cite{johnson91}
showed that $\bloneg$ is always weakly amenable. The operator space analogue 
of that result was established for Fourier algebras by the third author \cite{spronk},
and independently by the second author \cite{samei0}:  $\falg$ is always
operator weakly amenable. However, the failure of $\fal{\mathrm{SO}(3)}$ to be weakly amenable 
in the classical sense demonstrates once again that the operator space category was the right one to consider 
for weak amenability of $\falg$, just as it was for amenability of $\falg$. This time though, classical weak amenability
for $\falg$ did not require the group $G$ itself to be almost abelian. Indeed, as early as 1988, the 
first author showed $\falg$ to be weakly amenable for all discrete $G$. 
In \cite{johnson}, Johnson 
proved that $\falg$ is weakly amenable for any compact totally disconnected group.  This result was
extended to a wider class of locally compact groups
by the first author in \cite{forrest}:  if $G=A\times D$ where $A$ is abelian
and $D$ is totally disconnected, then $\falg$ is weakly amenable.  Together, these positive 
results pointed to the structure of the connected component $G_e$ of $G$ as the deciding factor 
in the question of the potential weak amenability of 
$\falg$. In fact, as early as 1996, and even prior to the appearance of \cite{forrest} in print, it 
had been conjectured that $\falg$ is weakly amenable precisely when $G_e$ is abelian. One direction of this 
conjecture was established by Forrest and Runde \cite{forrestr} who showed  that in 
fact $\falg$ is weakly amenable when $G_e$, the connected component of the identity
in $G$ is abelian. They also presented evidence for the converse when $G$ is a small invariant 
neighborhood group. 

The third author \cite{samei}
developed the theory of hyper-Tauberian commutative Banach algebras.
Though the hyper-Tauberian propety was origianally developed to study reflexivity of spaces of derivations
it implies weak amenability, and its theory
parallels the theory of weak amenability quite closely.  In particular, $\falg$
is always operator hyper-Tauberian, and is hyper-Tauberian when $G_e$ is 
abelian.  (Unfortunatly, the proof of the latter fact \cite[Theorem 21]{samei} contains a minor
error; this error is repaired in \cite[Theorem 3.7]{forrestss}.)
As is the case with weak amenability, it is thought likely that $\falg$ is hyper-Tauberian only when $G_e$
is abelian.

In the present article, we show that for compact groups, weak amenability of $\falg$
does indeed imply that $G_e$ is abelian. Using this result, we can establish the converse of 
the above mentioned weak amenability conjecture for all small inavriant neighborhood groups and for any maximally 
almost periodic group. Moreover, our result also shows that many non-abelian, non-compact connected
Lie groups do not have weakly amenable Fourier algebras. A similar statement regarding 
when $\falg$ is hyper-Tauberian necessarily follows as well.  

The paper is organised as follows.  We first deal with compact $G$.
In Section \ref{ssec:weakamen} we establish the main result of this article:
$\falg$ is weakly amenable only when $G_e$ is abelian.  Moreover, we 
show that if $\falg$ admits even a weakly amenable ideal, then $G_e$ must be abelian.
In Section \ref{ssec:synthesis} we show that there is a connection
between weak amenability of $\falg$ and spectral properties of the
antidiagonal $\check{\Del}$ for $\falgg$.  
We then examine the weak amenability, and amenability of $\falg$, in
relation to the same properties for $\fdelg$, in Section \ref{ssec:connection}.
The interesting algebras $\fdelg$ arose in the authors' previous 
paper \cite{forrestss}, and, as we further demonstrate here, hold much of the
structural information of the Fourier algebras of compact groups.
In the second part of the paper,
we use the results for compact groups to obtain similar
results for more general locally compact groups, including small invariant
neighbourhood groups and maximally almost periodic groups.  

The present study owes much of its motivation to work done by the authors 
in \cite{forrestss}, to an extend that we consider this article a successor
to that.  The second two authors were also motivated by techniques
they discovered whilst working on \cite{azimifardss}. Perhaps not surprisingly, operator 
spaces will again play a significant role in this investigation.

\subsection{Amenability, weak amenability and the hyper-Tauberian property}
We will be making non-trivial use
of operator spaces in Sections \ref{ssec:connection} and
\ref{ssec:locsynthesis}.  We will follow the definitions and notational
conventions of \cite{effrosrB}.  In particular, we use the ideas of
{\it completely bounded}, {\it completely contractive}, {\it completey isometric}
and {\it complete quotient} operators; {\it operator dual spaces}; and the {\it operator projective tensor
product} $\fV\what{\otimes}\fW$ of two complete operator spaces.
We use $\fX\otimes^\gam\fY$ to denote the {\it projective tensor product}
of two Banach spaces.  A {\it completely contractive Banach algebra}
is a complete operator space $\fA$, which admits an associative
product $m_\fA^0:\fA\otimes\fA\to\fA$ which extends to a complete
contraction $m_\fA:\fA\what{\otimes}\fA\to\fA$.  A complete operator
space and $\fA$-bimodule $\fV$ is a {\it completely contractive} $\fA$-bimodule
if the left and right actions $\fA\otimes\fV,\fV\otimes\fA\to\fV$
extend to complete contractions $\fA\what{\otimes}\fV,
\fV\what{\otimes}\fA\to\fV$.

Let $\fA$ be a (completely contractive) Banach algebra.  
Following \cite{johnsonM,ruan},
we say that $\fA$ is {\it (operator) amenable} if every
(completely) bounded derivation $D:\fA\to\fV^*$ -- i.e.\ linear
map for which $D(ab)=a\mult D(b)+D(a)\mult b$ -- for an
(completely) bounded dual bimodule $\fV$, is inner -- 
i.e.\ there is $f\iin\fV^*$
such that $Da=a\mult f-f\mult a$.  We say $\fA$ is {\it (operator)
weakly amenable} if every (completely) bounded derivation $D:\fA\to\fA^*$
is inner.  If $\fA$ is commutative, (operator) weak amenability
is equivalent to having no non-zero derivations $D:\fA\to\fV^*$
for every symmetric -- i.e.\ $a\mult v=v\mult a$ -- (completely)
bounded bimodule $\fV$ (see \cite[Theorem 1.5]{badecd} and 
\cite[Proposition 3.2]{forrestw}).

Let $\fA$ be a commutative (completely contractive) Banach algebra.
We suppose further that $\fA$ is a regular function algebra on a locally compact space $X$.
If $\vphi\in\fA^*$ we define
\[
\supp{\vphi}=\left\{x\in X:
\begin{matrix} \text{for every neighbourhood }U\oof x 
\text{ there is }f \\
\iin\fA\text{ such that }\supp{f}\subset U\aand \vphi(f)\not=0
\end{matrix}\right\}.
\]
Here $\supp{f}=\wbar{\{x\in X:f(x)\not=0\}}$.  An operator
$T:\fA\to\fA^*$ is called a {\it local map} if 
\[
\supp{Tu}\subset\supp{u}
\]
for every $u\iin\fA$.  We say $\fA$ is {\it (operator) hyper-Tauberian} on $X$ if 
every (completely bounded) bounded local map $T:\fA\to\fA^*$ is an 
$\fA$-module map.  This concept was developed by the third author
in \cite{samei}.  In that article, all results were stated for the case that
$X$ is the Gelfand spectrum of $\fA$; however, many of the proofs are valid
in the more general assumptions that $\fA$ is a regular function algebra on $X$.
For example, if $\fA$ is (operator) hyper-Tauberian then it is 
Tauberian and (operator) weakly amenable by \cite[Theorem 5]{samei} (resp.\ 
\cite[Theorem 26(ii)]{samei}).  
However, it is possible for $\fA$ to be
weakly amenable but not hyper-Tauberian \cite[Remark 24(ii)]{samei}.

Now we suppose that $\fA$ is a regular (completely contractive) Banach function algebra
on a locally compact space $X$.
If $E$ is a closed subset of $X$ we let
\begin{align*}
\ideal_\fA(E)&=\{u\in \fA:u(x)=0\text{ for every }x\iin E\} \\
\ideal_\fA^c(E)&=\{u\in\ideal_\fA(E):\supp{u}\text{ is compact}\},\aand \\
\ideal_\fA^0(E)&=\{u\in \fA:\supp{u}\text{ is compact and }\supp{u}\cap E=
\varnothing\} 
\end{align*}
We say that $E$ is a set of {\it spectral synthesis}, or is simply 
{\it spectral}, for $\fA$, if $\ideal_\fA(E)=\wbar{\ideal_\fA^0(E)}$.
We say $\fA$ is {\it Tauberian} on $X$ if the empty set $\varnothing$, 
qua subset of $X$, is a spectral for $\fA$.
We say $E$ is a set of {\it local synthesis} if $\ideal_\fA^c(E)
\subset \wbar{\ideal_\fA^0(E)}$.  Note that if $\fA$ is unital, and hence $X$ is compact,
local synthesis is the same as spectral synthesis.
We say that $E$ is {\it approximable} for $\fA$ if
$\ideal_\fA(E)$ admits a bounded approximate identity.
We say $E$ is {\it essential} for $\fA$ if $\ideal_\fA(E)$ is an essential
module over itself, i.e.\ $\wbar{\ideal_\fA(E)^2}=\ideal_\fA(E)$.

Let us summarise, for convenience of the reader, 
some of the strong connections which exist between
the amenability-like properties, listed above, and properties of ideals.

\begin{ideals}\label{lem:ideals}
Let $\fA$ be a regular (completely contractive) Banach function
algebra on a locally compact space $X$, such that 
$\fA\otimes^\gam\fA$ (resp.\ $\fA\what{\otimes}\fA$) is semi-simple.
Then the following conditions hold.

{\bf (i)} $\fA\otimes^\gam\fA$ (resp.\ $\fA\what{\otimes}\fA$)
is regular on $X\cross X$.

{\bf (ii)} The product map $m_\fA:\fA\otimes^\gam\fA\to\fA$
(resp.\ $m_\fA:\fA\what{\otimes}\fA\to\fA$) satisfies
$\ker m_\fA=\ideal_{\fA\otimes^\gam\fA}(\Del)$ (resp.\ 
$\ker m_\fA=\ideal_{\fA\what{\otimes}\fA}(\Del)$)
where $\Del=\{(x,x):x\in X\}$.

{\bf (iii)}  $\fA$ is (operator) amenable if and only if 
$\fA$ has a bounded approximate identity and
$\Del$ is approximable for $\fA\otimes^\gam\fA$ (resp.\ 
$\fA\what{\otimes}\fA$).

{\bf (iv)} $\fA$ is (operator) hyper-Tauberian on $X$ if and only if
$\fA$ is Tauberian and
$\Del$ is a set of local synthesis for $\fA\otimes^\gam\fA$ 
(resp.\ $\fA\what{\otimes}\fA$).

{\bf (v)} If $\fA$ has a bounded approximate identity, then
$\fA$ is (operator) weakly amenable if and only if $\wbar{\fA^2}=\fA$
and $\Del$ is essential for 
$\fA\otimes^\gam\fA$ (resp.\ $\fA\what{\otimes}\fA$).
\end{ideals}

\proof (i) This is \cite[Theorem 3]{tomiyama}.  The semisimplicity
of $\fA\otimes^\gam\fA$ (resp.\ $\fA\what{\otimes}\fA$)
is sufficient to guarantee that it is a function algebra on $X\cross X$.

(ii) This is immediate.  

(iii) This follows from (i) and a splitting result
of \cite{helemskii} (see also \cite[Theorem 3.10]{curtisl}).
The proof of \cite[Theorem 3.10]{curtisl}
can be adapted work in the completely contractive case, see
\cite[Theorem 3]{wood} or \cite[Lemma 2.1]{forrests}.

(iv) This is \cite[Theorem 6]{samei}.  It is noted in the proof of
\cite[Thoerem 3.10]{forrestss} how to adapt this to the operator
space setting.  We note that these proofs work for our general assumptions
so we do not need that $X$ is the spectrum of $\fA$.

(v) This follows fom (ii) and \cite[Theorem 3.2]{groenbaek} (see \cite[Section 2]{spronk}
for the completely bounded version). \endpf

\subsection{The Fourier algebra and some of its algebras}
The Fourier algebra $\falg$, for any locally compact group, was defined
in \cite{eymard}.  $\falg$ gains it operator space structure from being
the predual of a von Neumann algebra; see \cite[(3.10)]{eymard} and then
\cite[Section 3.2]{effrosrB}.
We note that $\falg\otimes^\gam\falg$ is naturally
isomorphic to $\falgg$ only when $G$ has an abelian subgroup of finite index, 
by \cite{losert}.  However, $\falg\what{\otimes}\falg\cong\falgg$
completely isometrically, via the natural identification; see 
\cite[Theorem 7.2.4]{effrosrB}.

If $G$ is a compact group we let $\Del=\{(s,s):s\in G\}$ denote the diagonal
subgroup of $G\cross G$, and let 
\[
\fal{G\cross G\!:\!\Del}=\{u\in\falgg:r\mult u=u\ffor r\iin G\}
\]
where $r\mult u(s,t)=u(srtr)$.  Then $\fal{G\cross G\!:\!\Del}$ is the closed 
subspace consists of $\falgg$ which
of functions constant on left cosets of $\Del$.  There is a homeomorphism
of the left coset space $G\cross G/\Del$ with $G$ given by 
$(s,t)\Del\mapsto st^{-1}$.  We let
\[
\fdelg=\{u\iin\fC(G):u(s)=w(s,e)\text{ for some }
w\iin\fal{G\cross G\!:\!\Del}\}
\]
The map
\[
\Gam:\falgg\to\fdelg,\quad \Gam w(s)=\int_G w(st,t)dt
\]
is surjective, and is injective on $\fal{G\cross G\!:\!\Del}$.
We endow $\fdelg$ with the norm and operator space structure
which make $\Gam$ a complete quotient map.  We note then that
\begin{equation}\label{eq:Nmap}
N:\fdelg\to\fal{G\cross G\!:\!\Del},\quad Nu(s,t)=u(st^{-1})
\end{equation}
is thus a complete isometry.  We note that
$\fdelg\what{\otimes}\fdelg\cong\fdelgg$ by \cite[Proposition 2.5]{forrestss}.
If we repeat the procedure above we obtain 
\[
\fdelbg{2}=\Gam\bigl(\fdelgg\bigr).
\]
We can do a similar construction with $\check{\Del}=\{(s,s^{-1}):s\in G\}$.
We let $G\cross G/\check{\Del}$ denote the set of equivalence classes
modulo the equivalence relation $(s',t')\sim (s,t)$ if and only if
$(s's^{-1},t^{-1}t')\in\check{\Del}$, so $G\cross G/\check{\Del}$, with
the quotient topology, is homeomorphic to $G$ via
$(s,t)\check{\Del}\mapsto st$.  We let
\[
\fal{G\cross G\!:\!\check{\Del}}=\{u\in\falgg:r\diamond u=u\ffor r\iin G\}
\]
where $r\diamond u(s,t)=u(sr,r^{-1}t)$. Similarly as above,
$ \fal{G\cross G\!:\!\check{\Del}}$ is a closed subalgebra of $\falgg$.
Let
\[
\fgamg=\{u\iin\fC(G):u(s)=w(s,e)\text{ for some }
w\iin\fal{G\cross G\!:\!\check{\Del}}\}
\]
Then the map
\[
\check{\Gam}:\falgg\to\fgamg, \quad \Gam w(s)=\int_G w(st,t^{-1})dt
\]
is surjective, and is injective on $\fal{G\cross G\!:\!\check{\Del}}$.
We endow $\fgamg$ with the norm and operator space structure
which make $\check{\Gam}$ a complete quotient map.  We also note that
\[
\check{N}:\fgamg\to\fal{G\cross G\!:\!\check{\Del}},\quad Nu(s,t)=u(st)
\]
is a complete isometry. 

We note the origional construction of $\fgamg$, from \cite{johnson}.
Let us first observe that $\falg$ has the approximation property, being 
metrically an $\ell^1$-direct sum of finite dimensional trace class
algebras \cite[(34.4)]{hewittrII}, and thus $\falg\otimes^\gam\falg$ is semi-simple with spectrum
$G\cross G$, by \cite[Theorem 4]{tomiyama}.  Thus $\falg\otimes^\gam\falg$
can be regarded as a subalgebra of $\falgg$, and we can restrict the maps
$\Gam$ and $\check{\Gam}$ to $\falg\otimes^\gam\falg$. We have
\[
\Gam\bigl(\falg\otimes^\gam\falg\bigr)=
\check{\Gam}\bigl(\falg\otimes^\gam\falg\bigr)=\fgamg
\]
as was observed in \cite[Section 4.2]{forrestss}.  Moreover,
each of $\Gam,\check{\Gam}:\falg\otimes^\gam\falg\to\fgamg$ are quotient
maps.  In \cite{johnson}, $\fgamg$ is used in a careful and clever way
to study the amenability properties of $\falg$.  It is curious that 
the projective, and operator projective products of $\falg$ produce the same result, here.

It is easily shown (see \cite[Section 2.1 \& Section 4.2]{forrestss})
that
\[
\fdelbg{2}\subset\fgamg\subset\fdelg\subset\falg.
\]
with equality holding for any pair, and hence all pairs, exactly
when $G$ has an open abelian subgroup.
Moreover, the inclusion map are all contractions.

We note that each of $\falg$ and $\fdelg$ have spectrum $G$ -- see
\cite{eymard} and \cite[Proposition 1.1]{forrestss} respectively -- and $\falg$ is Tauberian. 
Since $\falgg$ is regular on $G\cross G$,
it follows that each of $\fdelg$ and $\fgamg$ is regular on $G$; inductively
it follows that $\fdelbg{2}$ is regular on $G$.  For example, if $K$ and $L$ are disjoint
compact subsets of $G$, then $K^*$ and $L^*$ are disjoint compact subsets of $G\cross G$ where
$E^*=\{(s,t)\in G\cross G:st^{-1}\in E\}$ for $E\subset G$. Find $w\in\falgg$ such that
$w|_{K^*}=1$ and $w|_{L^*}=0$.  Then $u=\Gam w$ satisfies $u|_K=1$ and $u|_L=0$, showing
that $\fdelg$ is regular on $G$.
Though we suspect it to be true, we have not verified that either $\fgamg$ or $\fdelbg{2}$
have spectrum $G$.

\section{Compact groups}

\subsection{Weak amenability of Fourier algebras}\label{ssec:weakamen}
We let $G_e$ denote the connected component of the identity in $G$.
The following result, for compact groups, is the
the converse of \cite[Theorem 3.3]{forrestr}.

\begin{weakamen}\label{theo:weakamen}
Let $G$ be a compact group.  Then the following are equivalent:

{\bf (i)} $\falg$ is weakly amenable;

{\bf (ii)} $\falg$ is hyper-Tauberian;

{\bf (iii)} $G_e$ is abelian.
\end{weakamen}

\proof That (iii) implies (ii) is \cite[Theorem 21]{samei}.  That
(ii) implies (i) is \cite[Theorem 5]{samei}.  (We note that (iii) implies
(i) follows from \cite[Theorem 3.3]{forrestr}.)  Thus it suffices to
show that (i) implies (iii).

The restriction map $u\mapsto u|_{G_e}$ from $\falg$ to $\fal{G_e}$
is a surjection by \cite{herz}.  Thus, by \cite{dalesB}[Proposition 2.8.64],
if $\falg$ is weakly amenable, then so too must be $\fal{G_e}$.
Hence we assume that $G$ is connected, so $G=G_e$.

Suppose $G$ is connected and nonabelian.  Then by \cite[Theorem 6.5.6]{price}
there is a family $\{G_i\}_{i\in I}$ of compact connected Lie groups, at least one
of which is simple (in the sense of Lie groups) such that
\[
G\cong\left(\prod_{i\in I}G_i\right)/A
\]
where $A$ is a closed subgroup of the centre of $P=\prod_{i\in I}G_i$.
Let $G_{i_0}$ be simple.  Then, as shown in \cite{plymen}, there is 
a closed subgroup $K$ of $G_{i_0}$ such that
\[
\text{either }K\cong\mathrm{SU}(2)\oor K\cong\mathrm{SO}(3).
\]
Let $H=KA/A$, where $K$ is embedded into $P$ in the natural way.
By the second isomorphism theorem, $H\cong K/(K\cap A)$.  We have
$K\cap A\subset \z{K}$, the centre of $K$.  We note that
$\z{\mathrm{SO}(3)}$ is trivial, while $\z{\mathrm{SU}(2)}=\{-1,1\}$
and $\mathrm{SU}(2)/\{-1,1\}\cong\mathrm{SO}(3)$.  It follows that
$H$ is isomorphic to one of $\mathrm{SO}(3)$ or $\mathrm{SU}(2)$.
In either case, it follows from \cite{plymen} (via \cite[Corollary 7.3]{johnson})
that $\fal{H}$ is not weakly amenable.  Since $H$ is isomorphic to a closed 
subgroup of $G$, the restriction map gives a homomorphism from $\falg$ onto
$\fal{H}$.  Hence $\falg$ can not be weakly amenable either.  

Thus we conclude that for $\falg$ to be weakly amenable, $G_e$
must be abelian.  \endpf

\begin{idealweakamen}\label{theo:idealweakamen}
For a compact group $G$, if $\falg$ contains a weakly amenable
non-zero closed ideal, then $G_e$ is abelian.
\end{idealweakamen}

\proof  Let $\fI$ be a non-zero closed ideal, and 
\[
E=\hull(\fI)
=\{s\in G:u(s)=0\text{ for every }u\iin\fI\}.
\]
Suppose $\fI$ is weakly amenable.
Since left translation $u\mapsto s\con u$ ($s\con u(t)=u(s^{-1}t)$)
is an automorphism of $\falg$, for each $s\iin G$, each $s\con\fI$
is an ideal with $\hull(s\con\fI)=sE$.  For each finite
subset $F$ of $G$, the ideal $\fI_F=\wbar{\sum_{s\in F}s\con\fI}$
is weakly amenable.  Indeed, $\fI_F^*$ is a symmetric module
over each $s\con\fI$, so by \cite[Theorem 1.5]{badecd}
any bounded derivation $D:\fI_F\to\fI_F^*$
must satisfy $D|_{s\ast\fI}=0$ for each $s\iin F$.
Moreover, $\hull(\fI_F)=\cap_{s\in F}sE$.  By finite intersection
property we have either that

(a) $\cap_{s\in F}sE=\varnothing$ for some $F$, or

(b) $\cap_{s\in G}sE\not=\varnothing$.

\noindent In case (b), there is some  $t\in\cap_{s\in G}sE$, and thus 
$s\in Et^{-1}$ for every $s\iin G$, which implies that $E=G$.  
However, $\hull(\fI)=G$, only for $\fI=\{0\}$,
contradicting assumptions.  Hence case (a) holds, and for such $F$,
$\fI_F=\falg$, since $\falg$ is Tauberian.  But then $\falg$ itself
is weakly amenable, and $G_e$ is abelian by Theorem \ref{theo:weakamen}.
\endpf

\subsection{Spectral synthesis of the anti-diagonal}\label{ssec:synthesis}
Let $G$ be a compact group.

\begin{specsynth}\label{lem:specsynth}
Let $\theta:G\cross G\to G$ be given by $\theta(s,t)=st^{-1}$, and 
$\check{\theta}:G\cross G\to G$ be given by $\theta(s,t)=st$.
Then for a closed subset $E$ of $G$, the following are equivalent:

{\bf (i)} $E$ is spectral [approximable/essential] for $\fgamg$;

{\bf (ii)} $\check{\theta}^{-1}(E)$ is spectral [approximable/essential] for 
$\falgg$;

{\bf (iii)} $\theta^{-1}(E)$ is spectral [approximable/essential] for
$\falg\otimes^\gam\falg$;

{\bf (iii')} $\check{\theta}^{-1}(E)$ is spectral [approximable/essential] for
$\falg\otimes^\gam\falg$.
\end{specsynth}

\proof The equivalence of (i) and (iii) follows immediately from
\cite[Corollary 1.5]{forrestss}. 

The equivalence of (i) and (ii),
respectively of (i) and (iii'), follows similarly, but a few modifications are 
required.  In \cite[Theorem 1.4]{forrestss}
we implicitly used the group action 
of $G$ on $\falgg$, respectively on $\falg\otimes^\gam\falg$, 
given by $(r,w)\mapsto r\mult w$ (where $r\mult w(s,t)=w(sr,tr)$); see
\cite[Theorem 3.1]{spronkt} to see this action used more explicitly.
If the above group action is replaced 
by $(r,w)\mapsto r\diamond w$ (where $r\diamond w(s,t)=w(sr,r^{-1}t)$),
then it can be shown, similarly to \cite[Theorem 1.4]{forrestss} or
\cite[Theorem 3.1]{spronkt}, that
\[
\check{\Gamma}\ideal_{\fA(G\times G)}(\check{\theta}^{-1}(E))=
\ideal_{\falg}(E)
\quad\aand\quad
\langle\check{N}\ideal_{\falg}(E)\rangle=\ideal_{\fA(G\times G)}
(\check{\theta}^{-1}(E))
\]
where $\fA(G\cross G)$ is $\falgg$, respectively $\falg\otimes^\gam\falg$,
and $\langle \fS\rangle$ denotes the closed ideal generated by $\fS$.
The desired results then can be proved exactly as in 
\cite[Corollary 1.5]{forrestss}.  \endpf

In \cite{forrestr} it is shown that $\falg$ is amenable
exactly when $\check{\Del}=\{(s,s^{-1}):s\in G\}$ is an element of
the coset ring of $G\cross G$, exactly when $G$ admits an abelian subgroup
of finite index.  The next result is analogous to that.

\begin{antidiag}\label{theo:antidiag}
Let $G$ be a compact group.  Then the following are equivalent:

{\bf (i)} $\falg$ is weakly amenable;

{\bf (ii)} $\check{\Del}$ is spectral for $\falgg$;

{\bf (iii)} $\check{\Del}$ is essential for $\falgg$;

{\bf (iv)} $G_e$ is abelian.
\end{antidiag}

\proof The equivalence of (i) and (iv) is established in Theorem
\ref{theo:weakamen}.  That (ii) implies (iii) is clear.

If (iii) holds, then $\check{\theta}(\check{\Del})=\{e\}$ is essential
for $\fgamg$ by Lemma \ref{lem:specsynth},
(iii) implies (i).   It then follows, as in the proof of 
\cite[Theorem 7.2]{johnson}, that $\falg$ is weakly amenable.

If (iv) holds, then $\falg$ is hyper-Tauberian by \cite[Theorem 21]{samei}.
But then by \cite[Theorem 6]{samei} $\Del$ is a set of local synthesis
for $\falg\otimes^\gam\falg$, hence a set of synthesis as $G$ is compact.
Thus by Lemma \ref{lem:specsynth}, (iii) implies (i), $\{e\}$
is a set of synthesis for $\fgamg$.  But by Lemma \ref{lem:specsynth},
(i) implies (ii), we get that $\check{\Del}$ is spectral for
$\falgg$.  \endpf

We finish this section by giving alternative proofs to
some results from \cite{johnson}.  While our proof is ultimately no more 
efficient than the proofs given there, it places the results in a more general 
framework.

\begin{johnson}\label{prop:johnson}
Let $G$ be a compact group. Then

{\bf (i)}  \cite[Theorem 3.2]{johnson}
$\falg$ is amenable if and only if $\ideal_{\fgamg}\{e\}$ has a bounded
approximate identity; and

{\bf (ii)} \cite[Theorem 7.2]{johnson} 
$\falg$ is weakly amenable if and only if $\wbar{\ideal_{\fgamg}\{e\}^2}
=\ideal_{\fgamg}\{e\}$.
\end{johnson}

\proof These are immediate consequences of the equivalence of (i) and (iii) in
Lemma \ref{lem:specsynth}, and Lemma \ref{lem:ideals}, (iii) and (v),
respectively.  \endpf

\subsection{The connection between cohomology of $\falg$ and 
cohomology of $\fdelg$}\label{ssec:connection}
We illustrate, in the next two theorems, a curious connection that
exists between cohomology of $\falg$, operator cohomology
of $\fdelg$ and cohomology of $\fdelg$.

\begin{weakamenfdelg}\label{theo:weakamenfdelg}
Let $G$ be a compact group.  Then the following are equivalent:

{\bf (i)} $\fdelg$ is operator weakly amenable;

{\bf (i')} $\fdelg$ is weakly amenable;

{\bf (ii)} $\fdelg$ is operator hyper-Tauberian;

{\bf (ii')} $\fdelg$ is hyper-Tauberian;

{\bf (iii)} $G_e$ is abelian.
\end{weakamenfdelg}

\proof  That (i') implies (i), and (ii') implies (ii) are clear.
That (ii') implies (i') is from \cite[Theorem 5]{samei}, and the 
completely bounded analogue, that (ii) implies (i), is noted in 
\cite[Theorem 26]{samei}.
That (iii) implies (ii') is from \cite[Theorem 3.7]{forrestss}
and the identification $\fdelg=\fal{G\cross G\!:\! \Del}$.
Thus it remains to show that (i) implies (iii).

If $\fdelg$ is operator weakly amenable, then 
since $\fdelg\what{\otimes}\fdelg\cong\fdelgg$ by 
\cite[Proposition 2.5 (i)]{forrestss}, it follows Lemma \ref{lem:ideals}
that 
\begin{equation}\label{eq:delidealsquare}
\wbar{\ideal_{\fdelbg{2}}\{e\}^2}=\ideal_{\fdelbg{2}}\{e\}.
\end{equation}

We will show that
\begin{equation}\label{eq:gamidealsquare}
\wbar{\ideal_{\fgamg}\{e\}^2}=\ideal_{\fgamg}\{e\}.
\end{equation}
Indeed, let us first note that 
\begin{equation}\label{eq:density}
\ideal_{\fdelbg{2}}\{e\}\text{ is dense in }\ideal_{\fgamg}\{e\}.  
\end{equation}
Since $\fdelbg{2}$ is dense in $\fgamg$,
for any $u\iin \ideal_{\fgamg}\{e\}$ there is a sequence
$(u_n)\subset \fdelbg{2}$ which converges to $u$.  We note
that 
\[
|u_n(e)|=|u_n(e)-u(e)|\leq\unorm{u_n-u}\leq\fgamnorm{u_n-u}
\overset{n\to\infty}{\longrightarrow}0.
\]
Thus if $u_n'=u_n-u_n(e)1$, then $(u_n')\subset\ideal_{\fdelbg{2}}\{e\}$
with $\lim_{n\to\infty}\fgamnorm{u_n'-u}=0$.  Then
(\ref{eq:gamidealsquare}) follows form (\ref{eq:delidealsquare})
and (\ref{eq:density}).

By Proposition \ref{prop:johnson} (ii),
(\ref{eq:gamidealsquare}) implies that $\falg$ is weakly amenable.
Thus it follows from Theorem \ref{theo:weakamen} that
$G_e$ is abelian.  \endpf

\begin{amenfdelg}\label{theo:amenfdelg}
Let $G$ be a compact group.  Then the following are equivalent:

{\bf (i)} $\fdelg$ is operator amenable;

{\bf (i')} $\fdelg$ is amenable;

{\bf (ii)} $\falg$ is amenable;

{\bf (iii)} $G$ admits an open abelian subgroup.
\end{amenfdelg}

\proof Clearly (i') implies (i).  
The equivalence of (ii) and (iii) is \cite[Theorem 2.3]{forrestr}.
By \cite[Corollary 2.4]{forrestss}, if (iii) holds then
$\falg=\fdelg$, isomorphically, and hence (i') holds.
Thus it suffices to show that (i) implies (ii).

If (i) holds,  then from Lemma \ref{lem:ideals} (iii) we have
that $\ideal_{\fdel{G\times G}}(\Del)$ has a bounded 
approximate identity.  But then it follows from 
\cite[Corollary 1.5 (i)]{forrestss} that $\ideal_{\fdelbg{2}}\{e\}$
has a bounded approximate identity.  From (\ref{eq:density}) it follows that
$\ideal_{\fgamg}\{e\}$ has a bounded approximate identity.  Then by
Proposition \ref{prop:johnson} (i), $\falg$ is amenable. \endpf

Comparing Theorems \ref{theo:weakamenfdelg} and \ref{theo:amenfdelg},
it seems that the operator cohomology of $\fdelg$ is trivial
exactly when the same is true of the cohomology of $\fdelg$, which, in turn, 
is the same as the vanishing of the the cohomology of $\falg$.
This leads to a natural question.

\parbox{.8\linewidth}{\it Is the operator space structure on
$\fdelg$ maximal?}

\noindent A positive answer to the above question would explain
many of the implications if Theorems  \ref{theo:weakamenfdelg} and 
\ref{theo:amenfdelg}.  Though it seems unlikely that the operator space
structure on $\fdelg$ is maximal, in general, there is further curiosity
related to this, which we indicate below.

As shown in \cite[Proposition 1.5]{forrestr}, the map
$u\mapsto \check{u}:\falg\to\falg$ (where $\check{u}(s)=u(s^{-1})$)
is completely bounded only when $G$ has an abelian subgroup of finite
index.  By \cite[Theorem 4.5]{forrestw}, this is equivalent to
$\falg$ having the maximal operator space structure.  However,
this map is a complete isometry only when $G$ is abelian; see
\cite[Proposition 3.1]{ilies} and \cite[Proposition 3.4]{runde1}.

\begin{inversion}\label{prop:inversion}
For any compact group $G$, 
the map $u\mapsto\check{u}$ is a complete isometry on $\fdelg$.
\end{inversion}

\proof  The map $\sigma:G\cross G\to G\cross G$ given by
$\sig(s,t)=(t,s)$ is an isomorphism, hence
$J_\sig:\falgg\to\falgg$, $J_\sig w=w\comp\sig$, is a completely isometric
homomorphism; see, for example, \cite[Proposition 3.1]{ilies}.
It is clear that $J_\sig(\fal{G\cross G\!:\!\Del})=\fal{G\cross G\!:\!\Del}$.
Let $N$ be as in (\ref{eq:Nmap}), so $N$ has inverse
$M:\fal{G\cross G\!:\!\Del}\to\fdelg$ given by $Mw(s)=w(s,e)$, which is
thus also a complete isometry.  We note that
\[
[M J_\sig Nu](s)=[J_\sig Nu](s,e)=Nu(e,s)\overset{\dagger}{=}Nu(s^{-1},e)
=u(s^{-1})=\check{u}(s)
\]
where we used that $Nu\in\fal{G\cross G\!:\!\Del}$ at $\dagger$.
Hence $u\mapsto\check{u}$, which factors as three complete contractions,
is itself a complete contraction.
This map is obviously its own inverse, hence is a complete
isometry.  \endpf

\section{Some non-compact groups}

\subsection{A general weak amenability result}\label{ssec:genweakamen}

\begin{genweakamen}\label{theo:genweakamen}
If $G$ is a locally compact group for which $\falg$ admits a non-zero
weakly amenable ideal, then $G$ contains no non-abelian compact connected
subgroups.  In particular, this holds for $\fI=\falg$.
\end{genweakamen}

\proof Let $\fI$ be a non-zero ideal of $\falg$, and $K$ be a connected
compact subgroup of $G$.
Since $E=\hull(\fI)\not=G$, there is $s\iin G$ such that $sE\not\supset K$, and
$s\con\fI$ is also a weakly amenable ideal.
Since the restriction map $u\mapsto u|_K:\falg\to\fal{K}$ is surjective, by
\cite{herz}, $\fJ=\wbar{s\con\fI|_K}$ is a closed ideal of $\fal{K}$.  The choice
of $s$ ensures that $\hull(\fJ)\not=K$.  Since $\fJ$ is the closure of the image of
a weakly amenable algebra, $\fJ$ is weakly amenable by \cite[Proposition 2.8.64]{dalesB}.  
It then follows from
Theorem \ref{theo:idealweakamen} that $K$ must be abelian.  \endpf

The above result while very broad, misses still some key examples
such as $\mathrm{SL}_2(\Ree)$, the $ax+b$-group and the Heisenberg group.  
As we mentioned before, it is conjectured
that $\falg$ can be weakly amenable only when $G_e$ is abelian. It may well be the case, that 
$G_{e}$ is abelian if $\falg$ has even a non-zero weakly amenable closed ideal. 

\subsection{Local synthesis of the anti-diagonal}\label{ssec:locsynthesis}
We begin with a lemma, indicating a local property of local synthesis.

\begin{localproperty}\label{lem:localproperty}
Let $\fA$ be a regular semi-simple commutative Banach algebra with
spectrum $\Omega$, and $E$ be a closed subset of $\Omega$.
If each point in $E$ admits a closed neighbourhood $N$ such that
$N\cap E$ is of local synthesis, then $E$ itself is of local synthesis.
\end{localproperty}

\proof Let $u\in \ideal^c(E)$.  By assumption, there exists, for each $x\in E$,
a closed neighbourhood $N_x$ of $x$ such that
\begin{equation}\label{eq:localu}
u\in\wbar{\ideal^0(E_x)}\wwhere E_x=N_x\cap E
\end{equation}
as $\ideal^c(E_x)\supset\ideal^c(E)\ni u$.  By regularity of $\fA$,
there is $h_x\in\fA$ whose support is contained in the interior of $N_x$
and for which there is an open  neighbourhood $U_x$ of $x$ such that $h_x|_{U_x}=1$.
By (\ref{eq:localu}), given $\eps>0$, there is $u_x\in\ideal^0(E_x)$ such that
$\norm{u-u_x}<\eps/\norm{h_x}$.  Hence it follows that
$\norm{h_xu-h_xu_x}<\eps$.  We have that $h_xu_x\in\ideal^0(E)$ since
$h_xu_x$ vanishes on a neighbourhood of $E_x\cup(\Omega\setdif N_x)$,
and hence on a neighbourhood of $E$.  Since $\eps$ is arbitrary, we have
that $h_xu\in\wbar{\ideal^0(E)}$ and $h_xu_x|_{U_x}=u_x|_{U_x}$.  
In other words, ``$u$ is locally contained
in $\wbar{\ideal^0(E)}$''.  Then, a standard partition of unity argument
shows that $u\in\wbar{\ideal^0(E)}$.  \endpf

With the above lemma we can obtain The following generalisation of Theorem 
\ref{theo:antidiag}.

\begin{locantidiag}\label{theo:locantidiag}
Let $G$ be a locally compact group which admits an open subgroup
$H\cong A\cross K$, where $A$ is an abelian group and $K$ is compact.
Then the following are equivalent:
 
{\bf (i)} $\falg$ is weakly amenable;
 
{\bf (ii)} $\falg$ is hyper-Tauberian;
 
{\bf (iii)} $\check{\Del}=\{(s,s^{-1}):s\in G\}$ is a set of local
synthesis for $\falgg$;
 
{\bf (iv)} $\ideal'_{\fal{G\times G}}(\check{\Del})
=\wbar{\falgg\ideal_{\fal{G\times G}}(\check{\Del})}$ is operator
weakly amenable;

{\bf (v)} $G_e$ is abelian.
\end{locantidiag}

\proof That (v) implies (ii) is \cite[Theorem 21]{samei}, and that (ii) implies (i) is
\cite[Theorem 5]{samei}.  If (i) holds, then by Theorem \ref{theo:genweakamen},
$K_e$ is abelian, hence $G_e=H_e\cong A\cross K_e$ is abelian, so (v) holds.

Suppose (v) holds.  Then $K_e$ is abelian, and hence
by Theorem \ref{theo:antidiag} $\check{\Del}_K$ is a set of synthesis
for $\fal{K\cross K}$.  Since $\check{\Del}_A$ is a subgroup of $A\cross A$,
it is a set of synthesis for $\fal{A\cross A}$ by \cite[Theorem 2]{herz}.  
Hence it follows (the operator space
analogue of) \cite[Theorem 1.6]{kaniuth}, that $\check{\Del}_A\cross\check{\Del}_K$
is a set of synthesis for $\fal{A\cross A}\what{\otimes}\fal{K\cross K}
\cong\fal{A\cross A\cross K\cross K}$.  By applying a suitable automorphism
it follows that $\check{\Del}_H\cong\check{\Del}_{A\times K}$ is
a set of synthesis for $\fal{H\times H}\cong\fal{A\cross K\cross A\cross K}$.
Hence $\check{\Del}_H$ is a set of synthesis for $\falgg$, since $H$ is open in $G$.
If $(s,s^{-1})\in\check{\Del}_G$, then $s=th$ for some $t\iin G$ and $h\iin H$.
Hence $(s,s^{-1})\in(t,e)\check{\Del}_H(e,t^{-1})$.  Letting $N_s=(t,e)H\cross H(e,t^{-1})$
we have that $\check{\Del}_G\cap N_s=(t,e)\check{\Del}_H(e,t^{-1})$ is a set
of (local) synthesis.  Hence by Lemma \ref{lem:localproperty} we obtain (iii).

If (iii) holds, then 
\[
\ideal'_{\fal{G\times G}}(\check{\Del})=\wbar{\ideal^c_{\fal{G\times G}}(\check{\Del})}
=\ideal^0_{\fal{G\times G}}(\check{\Del})
\]
by regularity of $G$, and thus $\ideal'_{\fal{G\times G}}(\check{\Del})$ is essential.
Since $\falgg$ is operator weakly amenable by
\cite{spronk}, it follows \cite[Proposition 3.5]{forrestw} (which is the
completely bounded version of \cite[Corollary 1.3]{groenbaek0}) that (iv)
holds.

Finally, suppose (iv) holds.  Then 
$\ideal'_{\fal{G\times G}}(\check{\Del}_G)|_{K\times K}
=\ideal_{\fal{K\times K}}(\check{\Del}_K)$ is operator weakly
amenable.  Hence $\ideal_{\fal{K\times K}}(\check{\Del}_K)$ is essential
by \cite[Lemma 3.1]{forrestw}. Thus by Lemma \ref{lem:specsynth} (ii) with $E=\{e\}$, and
Proposition \ref{prop:johnson} (ii),
$\fal{K}$ is weakly amenable.  Hence $K_e$ is abelian and so
$G_e\cong A_e\cross K_e$ is abelian, and we obtain (v).  \endpf

We say that a locally compact group $G$ 
has {\it small invariant neighbourhoods}
if there is a basis for the identity of neighbourhoods invariant
under inner automorphisms.  We say $G$ is {\it maximally almost periodic}
if it continuously embeds into a compact group.

\begin{locantidiag1}\label{cor:locantidiag1}
If $G$ is either a small invariant neighbourhood group, or a
maximally almost periodic group, then the equivalent conditions
of Theorem \ref{theo:locantidiag} hold.
\end{locantidiag1}

\proof Any locally compact group $G$ admits an open almost connected subgroup
$H$ by \cite[(7.3) \& (7.7)]{hewittrI}.  If $G$ is either a small invariant 
neighbourhood group, or is maximally almost periodic, clearly $H$ has the
same property.  Then by \cite[Theorem 2.9]{grosserm} $H$ is isomorphic
to a semi-direct product $V\ltimes_\eta K$ where $V$ is a vector group
and $K$ is compact with $\eta(K)$ a finite
subgroup of automorphisms on $V$.  (See \cite[pages 1486-1487]{palmer}.)
Then $V\cross\ker\eta$ is an open subgroup of $H$, and hence
of $G$.  Thus we may appeal directly to Theorem \ref{theo:locantidiag}.
\endpf

\medskip
Address: {\sc Department of Pure Mathematics, University of Waterloo,
Waterloo, ON\quad N2L 3G1, Canada} 

\medskip
E-mail addresses:  {\tt beforres@uwaterloo.ca, esamei@uwaterloo.ca,  
\linebreak nspronk@uwaterloo.ca}


\begin{thebibliography}{10}

\bibitem{azimifardss}
A.~Azimifard, E.~Samei, and N.~Spronk.
\newblock Amenability properties of the centres of group algebras.
\newblock Preprint, see {\tt arXiv:0805.3685}, 2008.

\bibitem{badecd}
W.G. Bade, P.C. Curtis, and H.G. Dales.
\newblock Amenability and weak amenability for {B}eurling and {L}ipschitz
  algebras.
\newblock {\em Proc. London Math. Soc.}, 55:359--377, 1987.



\bibitem{curtisl}
P.C. Curtis and R.J. Loy.
\newblock The structure of amenable {B}anach algebras.
\newblock {\em J. London Math. Soc.}, 40:89--104, 1989.

\bibitem{dalesB}
H.G. Dales.
\newblock {\em Banach Algebras and Automatic Continuity}, volume~24 of {\em London Math. Soc., New
  Series}.
\newblock Claredon Press, Oxford Univ. Press, New York, 2000.

\bibitem{effrosrB}
E.G. Effros and Z.-J. Ruan.
\newblock {\em Operator Spaces}, volume~23 of {\em London Math. Soc., New
  Series}.
\newblock Claredon Press, Oxford Univ. Press, New York, 2000.

\bibitem{eymard}
P.~Eymard.
\newblock L'alg\`{e}bre de {F}ourier d'un groupe localement compact.
\newblock {\em Bull. Soc. Math. France}, 92:181--236, 1964.

\bibitem{forrest1}
Forrest, Brian.
\newblock Amenability and Derivations of the {F}ourier algebra.
\newblock {\em Proc. Amer. Math. Soc.}, 2:437--442, 1988.

\bibitem{forrest}
Forrest, Brian.
\newblock Fourier analysis on coset spaces.
\newblock {\em Rocky Mountain J. Math.}, 28:173--190, 1998.


\bibitem{forrestr}
B.E. Forrest and V.~Runde.
\newblock Amenability and weak amenability of the {F}ourier algebra.
\newblock {\em Math. Z.}, 250:731--744, 2005.

\bibitem{forrestss}
B.E. Forrest, E.~Samei, and N.~Spronk.
\newblock Convolutions on compact groups and {F}ourier algebras of coset
  spaces.
\newblock Preprint, see {\tt arXiv:0705.4277}, 2008.

\bibitem{forrests}
B.E. Forrest and N.~Spronk.
\newblock Best bounds for approximate identitiesin ideals of the {F}ouirier
  algebra vanishing on subgroups.
\newblock {\em Proc. Amer. Math. Soc.}, 134:111--116, 2005.

\bibitem{forrestw}
B.E. Forrest and P.J. Wood.
\newblock Cohomology and the operator space structure of the {F}ourier algebra
  and its second dual.
\newblock {\em Indiana Univ. Math. J.}, 50:1217--1240, 2001.

\bibitem{groenbaek}
N.~Groenbaek.
\newblock A characterization of weakly amenable {B}anach algebras.
\newblock {\em Studia Math.}, 94:149--162, 1989.

\bibitem{groenbaek0}
N.~Groenbaek.
\newblock Commutative {B}anach algebras, module derivations and semigroups.
\newblock {\em J. London Math. Soc. (2)}, 40:137--157, 1989.

\bibitem{grosserm}
S.~Grosser and M.~Moskowitz.
\newblock Compactness conditions in topological groups.
\newblock {\em J. Reine Angew. Math.}, 246:1--40, 1971.

\bibitem{herz}
C.S. Herz.
\newblock Harmonic synthesis for subgroups.
\newblock {\em Ann. Inst. Fourier, Grenoble}, 23(3):91--123, 1973.

\bibitem{hewittrI}
E.~Hewitt and K.A. Ross.
\newblock {\em Abstract Harmonic Analysis I}, volume 115 of {\em Grundlehern
  der mathemarischen Wissenschaften}.
\newblock Springer, New York, second edition, 1979.

\bibitem{hewittrII}
E.~Hewitt and K.A. Ross.
\newblock {\em Abstract Harmonic Analysis II}, volume 152 of {\em Grundlehern
  der mathemarischen Wissenschaften}.
\newblock Springer, New York, 1970.

\bibitem{ilies}
M.~Ilie and N.~Spronk.
\newblock Completely bounded homomorphisms of the {F}ourier algebras.
\newblock {\em J. Funct. Anal.}, 225:480--499, 2005.

\bibitem{johnsonM}
B.~E. Johnson.
\newblock {\em Cohomology in Banach algebras}, volume 127 of {\em Memoirs of
  the Amer. Math. Soc.}
\newblock 1972.

\bibitem{johnson91}
B.~E. Johnson.
\newblock Weak amenability of group algebras.
\newblock {\em Bull. London Math. Soc.}, 23:281--284, 1991.

\bibitem{johnson}
B.E. Johnson.
\newblock Non-amenability of the {F}ourier algebra of a compact group.
\newblock {\em J. London Math. Soc.}, 50:361--374, 1994.

\bibitem{kaniuth}
E.~Kaniuth.
\newblock Weak spectral synthesis for the projective tensor product of
  commutative banach algebras.
\newblock {\em Proc. Amer. Math. Soc.}, 132:2959--2967, 2004.

\bibitem{helemskii}
A.Ya. Khelemskii.
\newblock Flat {B}anach modules and amenable algebras.
\newblock {\em Trudy Moskov. Mat. Obshch.}, 47:179--218, 1984.

\bibitem{losert}
V.~Losert.
\newblock On tensor products of the {F}ourier algebras.
\newblock {\em Arch. Math.}, 43:370--372, 1984.


\bibitem{palmer}
T.W. Palmer.
\newblock {\em Banach algebras and the aeneral theory of $*$-algebras. Volume
  II. $*$-algebras}, volume~79 of {\em Encyclopedia of Mathematics and its
  applications}.
\newblock Cambridge University Press, Cambridge, 2001.

\bibitem{plymen}
R.J. Plymen.
\newblock {Fourier algebra of a compact Lie group}.
\newblock {\tt arXiv:math.FA/0104018}.

\bibitem{price}
J.~F. Price.
\newblock {\em Lie Groups and Compact Groups}, volume~25 of {\em London Math.
  Soc. Lec. Note Ser.}
\newblock Cambridge, 1977.

\bibitem{ruan}
Z.-J. Ruan.
\newblock The operator amenability of ${A(G)}$.
\newblock {\em Amer. J. Math.}, 117:1449--1474, 1995.

\bibitem{runde1}
V.~Runde.
\newblock The amenability constant of the fourier algebra.
\newblock {\em Proc. Amer. Math. Soc.}, 134:1473--1481, 2006.

\bibitem{samei0}
E.~Samei.
\newblock Bounded and completely bounded local derivations from certain
  commutative semisimple {B}anach algebras.
\newblock {\em Proc. Amer. Math. Soc.}, 133:229--238, 2005.

\bibitem{samei}
E.~Samei.
\newblock {Hyper-Tauberian algebras and weak amenability of
  Fig\`{a}-Talamanca-Herz algebras}.
\newblock {\em J. Funct. Anal.}, 231:195--220, 2006.

\bibitem{spronk}
N.~Spronk.
\newblock Operator weak amenability of the {F}ourier algebra.
\newblock {\em Proc. Amer. Math. Soc.}, 130:3609--3617, 2002.

\bibitem{spronkt}
N.~Spronk and L.~Turowska.
\newblock Spectral synthesis and operator synthesis for compact groups.
\newblock {\em J. London Math. Soc. (2)}, 66:361--376, 2002.

\bibitem{tomiyama}
J.~Tomiyama.
\newblock Tensor products of commutative {B}anach algebras.
\newblock {\em T\^{o}hoku Math. J.}, 12:143--154, 1960.

\bibitem{wood}
P.J. Wood.
\newblock Complemented idelas in the Fourier algebra of a locally compact
  group.
\newblock {\em Proc. Amer. Math. Soc.}, 128:445--451, 2000.

\end{thebibliography}
\end{document}